\noindent\centerline{\bf Erd\'elyi-Kober Fractional Integral Operators from a Statistical Perspective -IV}

\vskip.3cm\centerline{A.M. MATHAI}
\vskip.2cm\centerline{Centre for Mathematical Sciences,}
\vskip.1cm\centerline{Arunapuram P.O., Pala, Kerala-68674, India, and}
\vskip.1cm\centerline{Department of Mathematics and Statistics, McGill University,}
\vskip.1cm\centerline{Montreal, Quebec, Canada, H3A 2K6}
\vskip.2cm\centerline{and}
\vskip.2cm\centerline{H.J. HAUBOLD}
\vskip.1cm\centerline{Office for Outer Space Affairs, United Nations}
\vskip.1cm\centerline{P.O. Box 500, Vienna International Centre} 
\vskip.1cm\centerline{A - 1400 Vienna, Austria, and}
\vskip.2cm\centerline{Centre for Mathematical Sciences,}
\vskip.1cm\centerline{Arunapuram P.O., Pala, Kerala-68674, India}

\vskip.5cm\noindent{\bf Abstract}

\vskip.3cm In the preceding articles we considered fractional integral transforms involving one real scalar variable, one real matrix variable and real scalar multivariable case. In the present paper we consider the multivariable case when the arbitrary function is a real-valued scalar function of many $p\times p$ real matrix variables $X_1,...,X_k$. Extension of all standard fractional integral operators to the cases of many matrix variables is considered, along with interesting special cases and generalized matrix transforms.

\vskip.3cm\noindent{\bf 1.\hskip.3cm Introduction}

\vskip.3cm In all the preceding papers in this sequence the basic claim is that fractional integral operators are of two kinds, the first kind and the second kind. The first kind of operators belong to the class of Mellin convolution of a ratio and the second kind of operators belong to the class of Mellin convolution of a product. We will give the following formal definition of fractional integral operators of the first kind and second kind. The following standard notations will be used in this article. All the matrices appearing are $p\times p$ real symmetric, and further, positive definite (denoted by $X>O$) unless otherwise specified. $|X|$ denotes the determinant of the matrix $X$, ${\rm tr}(X) $ denotes the trace of $X$,

$$\eqalignno{{\rm d}X&=\prod_{i,j}\wedge{\rm d}x_{ij}\hbox{  for a general $X$ }\cr
&=\prod_{i\ge j}\wedge{\rm d}x_{ij}=\prod_{j\ge i}\wedge{\rm d}x_{ij}\hbox{ when $X$ is symmetric, $X=X'$}\cr}
$$that is, the wedge product of all differentials in $X$. $\int_{O<A<X<B}f(X){\rm d}X=\int_{O}^{B}f(X){\rm d}X$ will mean the integral over all positive definite $X$, such that $A>O,X>O,X-A>O,B>O,B-X>O$ where $A$ and $B$ are constant matrices. In this notation $O<X<I$ will indicate that all eigenvalues of $X$ are in the open interval $(0,1)$. All functions of matrix argument considered in this paper are real-valued scalar functions, whether the argument is one matrix or more matrices. The real matrix-variate gamma function will be denoted and defined as

$$\Gamma_p(\alpha)=\pi^{{p(p-1)}\over4}\Gamma(\alpha)\Gamma(\alpha-{1\over2})...\Gamma(\alpha-{{p-1}\over2}),
\Re(\alpha)>{{p-1}\over2}.\eqno(1.1)
$$Also, $X^{1\over2}$ will denote the positive definite square root of a real positive definite $p\times p$ matrix $X$. A prime will denote the transpose, that is, $X=X'$ means $X$ is symmetric. The following standard Jacobians will be used frequently in this article. For more details and for more Jacobians see Mathai (1997).

$$\eqalignno{Y=AXA',|A|\ne 0, X=X'&\Rightarrow {\rm d}Y=|A|^{p+1}{\rm d}X&(1.2)\cr
Y=X^{-1}&\Rightarrow {\rm d}X=\cases{|Y|^{-2p}{\rm d}Y\hbox{ for a general $X$}\cr
|Y|^{-(p+1)}{\rm d}Y\hbox{ for $X=X'$}\cr}&(1.3)\cr}
$$

\vskip.3cm\noindent{\bf Definition 1.1.}\hskip.3cm{\it Fractional integral operators of the first kind in one scalar or matrix variable case}\hskip.3cm A fractional integral of the first kind of order $\alpha$ and of one scalar or matrix variable is a Mellin convolution of a ratio with the first function $f_1(x_1)$ is of the form
$$f_1(x_1)={{\phi_1(x_1)(1-x_1)^{\alpha-1}}\over{\Gamma(\alpha)}},\Re(\alpha)>0\eqno(1.4)
$$in the real scalar variable case,
$$f_1(X_1)={{\phi(X_1)|I-X_1|^{\alpha-{{p+1}\over2}}}\over{\Gamma_p(\alpha)}},\Re(\alpha)>{{p-1}\over2}\eqno(1.5)
$$in the single matrix variable case,
where $\phi_1(x_1)$ or $\phi(X_1)$ is a specified function,  and

$f_2(x_2)=\phi_2(x_2)f(x_2)$ where $\phi_2(x_2)$ is a specified function and $f(x_2)$ is an arbitrary function, in the real scalar case and $f_2(X_2)=\phi_2(X_2)f(X_2)$ in the real matrix case where $X_2$ is a $p\times p$ real positive definite matrix, so that the fractional integral of the first kind of order $\alpha$ is given by the Mellin convolution formula for a ratio, namely
$$\eqalignno{g(u)&=\int_{v<u}\phi_1({{v}\over{u}})(1-{{v}\over{u}})^{\alpha-1}\phi_2(v){{v}\over{u^2}}f(v){\rm d}v&(1.6)\cr
\noalign{\hbox{in the scalar variable case, and for the matrix variable case}}
g(U)&=\int_{V<U}\phi_1(V^{1\over2}U^{-1}V^{1\over2}){{1}\over{\Gamma_p(\alpha)}}
|I-V^{1\over2}U^{-1}V^{1\over2}|^{\alpha-{{p+1}\over2}}\cr
&\times |V|^{{p+1}\over2}|U|^{-(p+1)}f(V){\rm d}V.&(1.7)\cr}
$$

\vskip.3cm\noindent{\bf Definition 1.2.}\hskip.3cm{\it Fractional integral operator of the second kind of order $\alpha$ for the one variable real scalar case or one matrix in the real case}\hskip.3cm Let $f_1$ and $f_2$ be as defined in Definition 1.1. Then the fractional integral operator of the second kind of order $\alpha$ is defined as the Mellin convolution of a product and given by

$$\eqalignno{g(u)&=\int_{v>u}{{1}\over{v}}\phi_1({{u}\over{v}}){{1}\over{\Gamma(\alpha)}}
(1-{{u}\over{v}})^{\alpha-1}\phi_2(v)f(v){\rm d}v&(1.8)\cr
\noalign{\hbox{for the real scalar case and for the real matrix case}}
g(U)&={{1}\over{\Gamma_p(\alpha)}}\int_{V>U}|V|^{-{{p+1}\over2}}\phi_1(V^{-{1\over2}}UV^{-{1\over2}})
|I-V^{-{1\over2}}UV^{-{1\over2}}|^{\alpha-{{p+1}\over2}}
\phi_2(V)f(V){\rm d}V.&(1.9)\cr}
$$

\vskip.3cm Note that if $f_1(x_1)$ and $f_2(x_2)$ are statistical densities of real positive scalar random variables $x_1$ and $x_2$ then $g(u)$ will represent the density of the ratio $u={{x_2}\over{x_1}}$. In the matrix case $g(U)$ will represent the density of $U$, where $X_1=V^{1\over2}U^{-1}V^{1\over2}$ and $X_2=V$. This is the statistical connection and later we will see that Kober operators are constant multiples of statistical densities.

\vskip.3cm\noindent{\bf Special cases}
\vskip.3cm\noindent{\bf Case (1)}:\hskip.3cm Let
$$\eqalignno{\phi_1(x)&=x_1^{\zeta-1},~\phi_2(x_2)=1\hbox{  then }\cr
g(u)&={{1}\over{\Gamma(\alpha)}}\int_v({{v}\over{u}})^{\zeta-1}(1-{{v}\over{u}})^{\alpha-1}{{v}\over{u^2}}f(v){\rm d}v\cr
&={{u^{-\zeta-\alpha}}\over{\Gamma(\alpha)}}\int_{v<u}(u-v)^{\alpha-1}v^{\zeta}f(v){\rm d}v\cr
&=\hbox{ Kober operator of the first kind for $\Re(\alpha)>0$}&(1.10)\cr}
$$In this case it is easily seen that ${{\Gamma(\zeta)}\over{\Gamma(\zeta+\alpha)}}g(u)$ is the density of $u={{x_2}\over{x_1}}$ where $x_1$ and $x_2$ are statistically independently distributed real scalar positive random variables.
\vskip.3cm\noindent{\bf Case (2)}:\hskip.3cm Let
$$\eqalignno{\phi_1(x_1)&=x_1^{\zeta-1},~\phi_2(x_2)=x_2^{-\zeta},~\zeta=-\alpha.\hbox{  Then}\cr
g(u)&={{1}\over{\Gamma(\alpha)}}\int_{v=a}^u(u-v)^{\alpha-1}f(v){\rm d}v={_aD_x^{-\alpha}}f(x)\cr
&={_aI_x^{\alpha}}f(x)\cr
&=\hbox{ Riemann-Liouville left sided fractional integral of order $\alpha$}.&(1.11)\cr}
$$If $a=0$ then we have a special case of the Riemann-Liouville left sided integral operator. When $a\to -\infty$ then it gives ${_{-\infty}W_x^{-\alpha}}f(x)= $ left sided Weyl operator.
\vskip.3cm\noindent{\bf Case (3):}\hskip.3cm
$$\phi_1(x_1)=x_1^{\zeta-1}{_2F_1}(\alpha+\beta,-\gamma;\alpha;(1-{{v}\over{u}}))
$$then $g(u)$ gives the Saigo operator of the first kind.
\vskip.3cm\noindent{\bf Case (4):}\hskip.3cm

$$\eqalignno{\phi_1(x_1)&=x_1^{\zeta},\phi_2(x_2)=1\hbox{  then  }\cr
g(u)&={{1}\over{\Gamma(\alpha)}}\int_v{{1}\over{v}}({{u}\over{v}})^{\zeta}(1-{{u}\over{v}})^{\alpha-1}f(v){\rm d}v\cr
&={{u^{\zeta}}\over{\Gamma(\alpha)}}\int_{v>u}v^{-\zeta-\alpha}(v-u)^{\alpha-1}f(v){\rm d}v&(1.12)\cr
&=K_{u}^{\zeta,\alpha}f(u)=\hbox{ Kober operator of the second kind}&(1.13)\cr}
$$This can be interpreted as a statistical density. In fact, ${{\Gamma(\zeta+1+\alpha)}\over{\Gamma(\zeta+1)}}g(u)$ is a statistical density of the product $u=x_1x_2$ where $x_1$ and $x_2$ are two positive statistically independent real scalar random variables with densities $f_1(x_1)$ and $f_2(x_2)$ respectively where $f_1$ is a type-1 beta density with the parameters $(\zeta+1,\alpha)$.

\vskip.3cm\noindent{\bf Case (5):}\hskip.3cm
$$\eqalignno{\phi(x_1)&=1,\phi_2(x_2)=x_2^{\alpha}\hbox{  then  }\cr
g(u)&={{1}\over{\Gamma(\alpha)}}\int_v{{1}\over{v}}(1-{{u}\over{v}})^{\alpha-1}v^{\alpha}f(v){\rm d}v={{1}\over{\Gamma(\alpha)}}\int_{v>u}(v-u)^{\alpha-1}f(v){\rm d}v\cr
&={_xD_{\infty}^{-\alpha}}f={_xW_{\infty}^{-\alpha}}\cr
&=\hbox{right sided Weyl fractional integral of order $\alpha$}.&(1.14)\cr}
$$Note that we can have the  matrix-variate cases also corresponding to the cases (1) to (5) above. Since they are straight forward the details are not given here.

\vskip.3cm\noindent{\bf Definition 1.3.}\hskip.3cm{\it Fractional integral operator of the first kind of orders $\alpha_1,...,\alpha_k$ in the multivariable case}\hskip.3cm Let $f_1(x_1,...,x_k)$ and $f_2(v_1,..,v_k)$ be real-valued scalar functions of the scalar variables $x_1,...,x_k$ and $v_1,...,v_k$ respectively, where $f_1$ is of the form

$$f_1=\phi_1(x_1,...,x_k)\{\prod_{j=1}^k{{(1-x_j)^{\alpha_j-1}}\over{\Gamma(\alpha_j)}}\}.
 $$Then the fractional integral operator of the first kind of orders $\alpha_1,...,\alpha_k$ in the multivariable scalar case is given by
 $$\eqalignno{g(u_1,...,u_k)&=\{\prod_{j=1}^k{{1}\over{\Gamma(\alpha_j)}}\int_{v_j<u_j}
 (1-{{v_j}\over{u_j}})^{\alpha_j-1}{{v_j}\over{u_j^2}}\}\cr
 &\times f(v_1,...,v_k){\rm d}V,{\rm d}V={\rm d}v_1\wedge...\wedge{\rm d}v_k.&(1.15)\cr}
 $$In the corresponding many matrix variable case $f_1$ is given by
 $$\eqalignno{f_1(U_1,...,U_k)&=\phi_1(V_1^{1\over2}U_1^{-1}V_1^{1\over2},...,V_k^{1\over2}U_k^{-1}V_k^{1\over2})\cr
 &\times\{\prod_{j=1}^k{{1}\over{\Gamma_p(\alpha_j)}}\int_{V_j<U_j}|I-V_j^{1\over2}U_j^{-1}V_j^{1\over2}|^{\alpha_j-1}
 |V_j|^{{p+1}\over2}|U|^{-(p+1)}\}\cr
 &\times f_2(V_1,...,V_k){\rm d}V,{\rm d}V={\rm d}V_1\wedge...\wedge{\rm d}V_k.&(1.16)\cr}
 $$

 \vskip.3cm\noindent{\bf Definition 1.4.}\hskip.3cm{\it Fractional integral operator of the second kind of orders $\alpha_1,...,\alpha_k$ in the multivariable case}\hskip.3cm Let $f_1(x_1,...,x_k)$ and $f_2(v_1,..,v_k)$ be real-valued scalar functions of the scalar variables $x_1,...,x_k$ and $v_1,...,v_k$ respectively, where $f_1$ is of the form
$$f_1=\phi_1(x_1,...,x_k)\{\prod_{j=1}^k{{(1-x_j)^{\alpha_j-1}}\over{\Gamma(\alpha_j)}}\}.
 $$Then the fractional integral operator of the second kind of orders $\alpha_1,...,\alpha_k$ in the multivariable scalar case is given by
 $$\eqalignno{g(u_1,...,u_k)&=\{\prod_{j=1}^k{{1}\over{\Gamma(\alpha_j)}}\int_{v_j>u_j}
 (1-{{u_j}\over{v_j}})^{\alpha_j-1}{{1}\over{v_j}}\}\cr
 &\times \phi_1({{u_1}\over{v_1}},...,{{u_k}\over{v_k}})f(v_1,...,v_k){\rm d}V,{\rm d}V={\rm d}v_1\wedge...\wedge{\rm d}v_k.&(1.17)\cr}
 $$In the corresponding many matrix variable case $g$ is given by

 $$\eqalignno{g(U_1,...,U_k)&=\int_{V_1>U_1}...\int_{V_k>U_k}\phi_1(V_1^{-{1\over2}}U_1V_1^{-{1\over2}},...,
 V_k^{-{1\over2}}U_kV_k^{-{1\over2}})\cr
 &\times\{\prod_{j=1}^k{{1}\over{\Gamma(\alpha_j)}}|I-V_j^{-{1\over2}}U_jV_j^{-{1\over2}}|^{\alpha_j-1}
 |V_j|^{-{{p+1}\over2}}\}\cr
 &\times f(V_1,...,V_k){\rm d}V,{\rm d}V={\rm d}V_1\wedge...\wedge{\rm d}V_k.&(1.18)\cr}
 $$
 \vskip.3cm We will give the following formal definitions for fractional derivatives. Fractional derivatives will be defined as the following fractional integrals. To this end let us have the following notations for the fractional integrals. All the following integrals are of order $\alpha$ in the single variable case and of orders $\alpha_1,...,\alpha_k$ in the multivariable case. Let

 $$\eqalignno{IF_x^{\alpha}&=\hbox{ fractional integral of the first kind, one variable case, eq: (1.6)}\cr
 IS_x^{\alpha}&=\hbox{ fractional integral of the second kind, one variable case, equation, eq: (1.8)}\cr
 IF_X^{\alpha}&=\hbox{ fractional integral of the first kind, one matrix variable case, eq: (1.7)}\cr
 IS_X^{\alpha}&=\hbox{ fractional integral of the second kind, one matrix variable case, eq: (1.9)}\cr}
 $$

 $$\eqalignno{IF_{x_1,...,x_k}^{\alpha_1,...,\alpha_k}&=\hbox{fractional integral of the first kind, eq: (1.15)}\cr
 IS_{x_1,...,x_k}^{\alpha_1,...,\alpha_k}&=\hbox{fractional integral of the second kind, eq: (1.17)}\cr
 IF_{X_1,...,X_k}^{\alpha_1,...,\alpha_k}&=\hbox{ fractional integral of the first kind, eq: (1.16)}\cr
 IS_{X_1,...,X_K}^{\alpha_1,...,\alpha_k}&=\hbox{ fractional integral of the second kind, eq: (1.18)}&(1.19)\cr}
 $$For all the following definitions $m=[\alpha]+1$ where $[\alpha]$ denotes the integer part of $\Re(\alpha)>0$ and $m_j=[\alpha_j]+1$ where $[\alpha_j]$ denotes the integer part of $\Re(\alpha_j)>0,~s j=1,...,k.$ Fractional integral can be written as an antiderivative and then, for example, the fractional integral $IF_x^{\alpha}$ can be written as an antiderivative as $IF_x^{\alpha}=DF_x^{-\alpha}$. Let $D_x^m$ denote the integer order derivative with respect to the real scalar variable $x$ and of order $m, m=0,1,2,...$. Then symbolically we can write
 $$D_x^{m}[IF_x^{m-\alpha}]=D_x^m[DF_x^{-(m-\alpha)}]=DF_x^{\alpha}\eqno(1.20)
 $$the corresponding fractional derivative of order $\alpha$. We can take (1.20) to define fractional derivative, where let $m=[\alpha]+1$, with $[\alpha]$ denoting the integer part of $\Re(\alpha)$. Let $D_j^{m_r}={{\partial^{m_r}}\over{\partial x_j}}$ denote the partial derivative of order $m_r$ with respect to $x_j$. Then, for example, we can define fractional derivatives of orders $\alpha_1,...,\alpha_k$ in the scalar multivariable case  as
 $$\eqalignno{D_1^{m_1}...D_k^{m_k}IF_{x_1,...,x_k}^{m_1-\alpha_1,...,m_k-\alpha_k}&=D_1^{m_1}...D_k^{m_k}
 DF_{x_1,...,x_k}^{-(m_1-\alpha_1),...,-(m_k-\alpha_k)}\cr
 &=DF_{x_1,...,x_k}^{\alpha_1,...,\alpha_k}.&(1.21)\cr}
 $$In the matrix-variate case we will introduce the following operator. Consider the matrix of partial derivatives, $${{\partial}\over{\partial X}}=({{\partial}\over{\partial x_{ij}}}),~ |{{\partial}\over{\partial X}}|\hbox{ determinant of ${{\partial}\over{\partial X}}$}\eqno(1.22)
 $$and the corresponding matrices and determinants of partial differential operators. Then in the single matrix-variate case we can define the fractional derivative of order $\alpha$ as follows:
 $$|{{\partial}\over{\partial X}}|^mIF_{X}^{m-\alpha}=|{{\partial}\over{\partial X}}|^{m}DF_X^{-(m-\alpha)}=DF_X^{\alpha}.\eqno(1.23)
 $$If there are several matrices then we consider the determiants of the matrices of partial differential operators as the differential operators and, for example, we can define
 $$\eqalignno{|{{\partial}\over{\partial X_1}}|^{m_1}...|{{\partial}\over{\partial X_k}}|^{m_k}&IF_{X_1,...,X_k}^{m_1-\alpha_1,...,m_k-\alpha_k}\cr
 &=|{{\partial}\over{\partial X_1}}|^{m_1}...|{{\partial}\over{\partial X_k}}|^{m_k}DF_{X_1,...,X_k}^{-(m_1-\alpha_1),...,-(m_k-\alpha_k)}\cr
 &=DF_{X_1,...,X_k}^{\alpha_1,...,\alpha_k}&(1.24)\cr}
 $$where $m_j=[\alpha_j]+1,j=1,...,k$. If any $\alpha_j$ is a positive integer then the corresponding derivative is a $(m_j-1)$th order total partial derivative and not a fractional one.

 \vskip.3cm\noindent{\bf Definition 1.5.}\hskip.3cm{\it Fractional derivative}\hskip.3cm Let the ``$I$'' in (1.19) be replaced by ``$D$'' to denote the corresponding fractional derivative. Then the fractional derivative of order $\alpha$ for the single scalar variable case will be defined as in (1.20) with the corresponding definitions for all other single variable fractional integrals in (1.19). The fractional derivative of order $\alpha$ in the single matrix case will be defined as in (1.23) with the corresponding changes for all the ssingle matrix variable cases in (1.19). Fractional derivatives of orders $\alpha_1,...,\alpha_k$ in the multivariable case will be defined as in (1.21) with the corresponding changes for all cases listed in (1.19), and in the many matrix variable case will be defined as in (1.24) with the corresponding changes for all the cases in (1.19).

 \vskip.3cm\noindent{\bf 2.\hskip.3cm Fractional Operators in the Many Matrix-variate Cases}

 \vskip.3cm Let $X_1,...,X_k$ and $V_1,...,V_k$ be two sequences of $p\times p$ matrix random variables where between the sets the two sets are statistically independently distributed. Further, let $X_1,...,X_k$ be mutually independently distributed type-1 real matrix-variate beta random variables with the parameters $(\zeta_j+{{p+1}\over2},\alpha_j),j=1,...,k$. That is, $X_j$ has the density

 $$\eqalignno{f_j(X_j)&={{\Gamma_p(\alpha_j+\zeta_j+{{p+1}\over2})}\over{\Gamma_p(\alpha_j)\Gamma_p(\zeta_j+{{p+1}\over2})}}
 |X_j|^{\zeta_j}
 |I-X_j|^{\alpha_j-{{p+1}\over2}}&(2.1)\cr}
 $$for $O<X_j<I$,  $\Re(\alpha_j)>{{p-1}\over2},\Re(\zeta_j)>-1$ and $f_j(X_j)=0$, $j=1,...,k$ elsewhere. Consider the transformation $V_j^{-{1\over2}}U_jV_j^{-{1\over2}},j=1,...,k$ then the Jacobian is $|V_1|^{-{{p+1}\over2}}...|V_k|^{-{{p+1}\over2}}$. Substituting in (2.1) the joint density of $U_1,...,U_k$, denoted by $g(U_1,...,U_K)$, is given by

 $$\eqalignno{g(U_1,...,U_k)&=\{\prod_{j=1}^k{{\Gamma_p(\alpha_j+\zeta_j+{{p+1}\over2})}
 \over{\Gamma_p(\alpha_j)\Gamma_p(\zeta_j+{{p+1}\over2})}}\cr
 \int_{V_j>U_j}|V_j|^{-{{p+1}\over2}}&|V_j^{-{1\over2}}U_jV_j^{-{1\over2}}|^{\zeta_j}
 |I-V^{-{1\over2}}U_jV_j^{-{1\over2}}|^{\alpha_j-{{p+1}\over2}}\}\cr
 &\times f(V_1,...,V_k){\rm d}V\cr
 &=\{\prod_{j=1}^k{{\Gamma_p(\alpha_j+\zeta_j+{{p+1}\over2})}\over{\Gamma_p(\alpha_j)\Gamma_p(\zeta_j+{{p+1}\over2})}}
 |U_j|^{\zeta_j}\cr
 &\times\int_{V_j>U_j}|V_j|^{-\zeta_j-\alpha_j}|V_j-U_j|^{\alpha_j-{{p+1}\over2}}\}f(V_1,...,V_k){\rm d}V.&(2.2)\cr}
 $$Hence we will define Kober operator of the second kind and of orders $(\alpha_1,...,\alpha_k)$ for the many matrix-variate case, and denoted as follows:

 $$\eqalignno{K_{U_1,...,U_k}^{(\zeta_j,\alpha_j),j=1,...,k}f(U_1,...,U_k)&=\{\prod_{j=1}^k
 {{|U_j|^{\zeta_j}}\over{\Gamma_p(\alpha_j)}}\int_{V_j>U_j}|V_j|^{-\zeta_j-\alpha_j}|V_j-U_j|^{\alpha_j-{{p+1}\over2}}\}\cr
 &\times f(V_1,...,V_k){\rm d}V.&(2.3)\cr}
 $$Therefore this Kober operator is a constant times a statistical density function, namely,
 $$\{{{\Gamma_p(\zeta_j+{{p+1}\over{2}})}\over{\Gamma_p(\alpha_j+\zeta_j+{{p+1}\over2})}}\}g(U_1,...,U_k)
 =K_{U_1,...,U_k}^{(\zeta_j,\alpha_j),j=1,...,k}f(U_1,...,U_k).\eqno(2.4)
 $$Now, let us consider $f_1$ having a joint density of the matrix variables $X_1,...,X_k$ and $f_2$ is a joint density of the matrix variables $V_1,...,V_k$ where the two sets are independently distributed. Then we can have several interesting results where the Kober operator of (2.4) will become constant multiples of statistical densities coming from various considerations.

 \vskip.3cm\noindent{\bf Theorem 2.1.}\hskip.3cm{\it Let the two sets $X_1,...,X_k$ and $V_1,...,V_k$ of matrices be independently distributed. Further, let $X_1,...,X_k$ have a joint type-1 Dirichlet density with the parameters $(\zeta_j+{{p+1}\over2},\alpha_j)$, $j=1,...,k$. Consider the transformation
 $$\eqalignno{X_1&=Y_1\cr
 X_2&=(I-Y_1)^{1\over2}Y_2(I-Y_1)^{1\over2}\cr
 X_j&=(I-Y_{j-1})^{1\over2}...(I-Y_1)^{1\over2}Y_j(I-Y_1)^{1\over2}...(I-Y_{j-1})^{1\over2},j=2,...,k.&(2.5)\cr
 \noalign{\hbox{Or}}
 Y_j&=(I-X_1-...-X_{j-1})^{-{1\over2}}X_j(I-X_1-...-X_{j-1})^{-{1\over2}},j=2,...,k, Y_1=X_1.&(2.6)\cr}
 $$Consider the transformation
 $$U_j=V_j^{1\over2}Y_jV_j^{1\over2},~Y_j=V_j^{-{1\over2}}U_jV_j^{-{1\over2}},j=1,...,k.\eqno(2.7)
 $$Then the joint density of $U_1,...,U_k$ is constant times the generalized Kober operator of the second kind defined in (2.4).}

 \vskip.3cm\noindent{\bf Proof:}\hskip.3cm Under the transformation in (2.5) or (2.6) the Jacobian is
 $$J=|I-Y_1|^{(k-1)({{p+1}\over2})}|I-Y_2|^{(k-2)({{p+1}\over2})}...|I-Y_k|^{{p+1}\over2}\eqno(2.8)
 $$and that $Y_1,...,Y_k$ are independently distributed as type-1 real matrix-variate beta random variables with the parameters $(\zeta_j+{{p+1}\over2},\beta_j),j=1,...,k$ where
 $$\beta_j=\zeta_{j+1}+\zeta_{j+2}+...+\zeta_k+(k-j)\eqno(2.9)
 $$see, for example, Mathai (1997). Now, it is equivalent to the situation in (2.3) and (2.4) with $Y_j$'s standing in place of the independently distributed $X_j$'s and hence from (2.4) we have the following result:

 $$\{\prod_{j=1}^k{{\Gamma_p(\zeta_j+{{p+1}\over2})}\over{\Gamma_p(\beta_j+\zeta_j+{{p+1}\over2})}}\}
 g(U_1,...,U_k)=K_{U_1,...,U_k}^{(\zeta_j,\beta_j),j=1,...,k}f(U_1,...,U_k)\eqno(2.10)
 $$where $\beta_j$ is given in (2.9). Hence the result.
 \vskip.3cm We can consider several generalized models belonging to the family of generalized type-1 Dirichlet family in the many matrix-variate cases. In all such situations we can derive the generalized Kober operator of the second kind in many matrices. We will take one such generalization here and obtain a theorem. Let $f_1(X_1,...,X_k)$ of the form

 $$\eqalignno{f_1(X_1,...,X_k)&=C~|X_1|^{\zeta_1}|I-X_1|^{\beta_1}|X_2|^{\zeta_2}|I-X_1-X_2|^{\beta_2}...\cr
 &\times|X_k|^{\zeta_k}|I-X_1-...-X_{k_1}|^{\zeta_k}|I-X_1-...-X_k|^{\beta_k-{{p+1}\over2}}&(2.11)\cr}
 $$where $O<X_1+...+X_j<I,j=1,...,k, \Re(\zeta_j)>-1,j=1,...,k$, $C$ is the normalizing constant. Other conditions on the parameters will be given later. In this connection we can establish the following theorem.

 \vskip.3cm\noindent{\bf Theorem 2.2.}\hskip.3cm{\it Let $X_1,...,X_k$ have a joint density as in (2.11). Consider the transformation as in (2.5) and (2.6) with the $U_j$'s  and $V_j$'s defined as in (2.7). Let the joint density of $U_1,...,U_k$ be denoted as $g(U_1,...,U_k)$. Let

 $$\delta_j=\zeta_{j+1}+...+\zeta_k+\beta_j+\beta_{j+1}+...+\beta_k,j=1,...,k.\eqno(2.12)
 $$Then
 $$\{\prod_{j=1}^k{{\Gamma_p(\zeta_j+{{p+1}\over2})}\over{\Gamma_p(\zeta_j+{{p+1}\over2}+\delta_j)}}\}
 g(U_1,...,U_k)=K_{U_1,...,U_k}^{(\zeta_j,\beta_j),j=1,..,k}f(U_1,...,U_k)\eqno(2.13)
 $$where $\delta_j$ is defined in (2.12) and the Kober operator in (2.4).}

 \vskip.3cm\noindent{\bf Proof:}\hskip.3cm We can show that under the transformation in (2.5) or (2.6) the $Y_j$'s are independently distributed as real matrix-variate type-1 beta random variables with the parameters $(\zeta_j+{{p+1}\over2},\delta_j),j=1,...,k$. Now the result follows from the procedure of the proof in Theorem 2.1. For Kober operators in the one variable case see Mathai and Haubold (2008).

 \vskip.3cm\noindent{\bf Note 2.1.}\hskip.3cm Special cases connecting to Riemann-Liouville operator, Weyl operator and Saigo operator, corresponding to the ones in section 1.1 for Kober operator of the second kind, can be obtained in a parallel manner and hence they will not be repeated here.

 \vskip.3cm\noindent{\bf 3.\hskip.3cm Fractional Integral Operators of the First Kind in the Case of Many Matrix Variables}

 \vskip.3cm Let $f_1(X_1,...,X_k)$ be a function of many $p\times p$ matrices and $f_2(V_1,...,V_k)$ be another function of another sequence of $p\times p$ matrices $V_1,...,V_k$. Let $f_1$ be of the form

 $$\{\prod_{j=1}^k{{\Gamma_p(\zeta_j+\alpha_j)}\over{\Gamma_p(\alpha_j)\Gamma_p(\zeta_j)}}
 |X_j|^{\zeta_j-{{p+1}\over2}}|I-X_j|^{\alpha_j-{{p+1}\over2}}\}\eqno(3.1)
 $$for $O<X_j<I,\Re(\alpha_j)>{{p-1}\over2},\Re(\zeta_j)>{{p-1}\over2},j=1,...,k.$ Let $X_j=V^{1\over2}U_j^{-1}V^{1\over2},j=1,...,k$. Let

$$\eqalignno{\{\prod_{j=1}^k{{\Gamma_p(\zeta_j)}\over{\Gamma_p(\zeta_j+\alpha_j)}}\}g(U_1,...,U_k)&
=\{\prod_{j=1}^k{{1}\over{\Gamma_p(\alpha_j)}}\int_{V_j}|V_j^{1\over2}U_j^{-1}V_j^{1\over2}|^{\zeta_j-{{p+1}\over2}}\cr
&\times |I-V_j^{1\over2}U_j^{-1}V_j^{1\over2}|^{\alpha_j-{{p+1}\over2}}|V_j|^{{p+1}\over2}|U_j|^{-(p+1)}\}\cr
&\times f(V_1,...,V_k){\rm d}V&(3.2)\cr
&=I_{U_1,...,U_k}^{(\zeta_j,\alpha_j),j=1,...,k}f(U_1,...,U_k).&(3.3)\cr}
$$Then (3.3) will be taken as the definition of fractional integral operator of the first kind of orders $\alpha_1,...,\alpha_k$ in the many matrix variable case.

\vskip.2cm Simplifying (3.2), we have a definition for generalized Kober operator of the first kind in many matrix variables case. We note that

$$\eqalignno{I_{U_1,...,U_k}^{(\zeta_j,\alpha_j),j=1,...,k}f(U_1,...,U_k)&
=\{\prod_{j=1}^k{{|U_j|^{-\zeta_j-\alpha_j}}\over{\Gamma_p(\alpha_j)}}\int_{V_j<U_j}|V_j|^{\zeta_j}
|U_j-V_j|^{\alpha_j-{{p+1}\over2}}\}\cr
&\times f(V_1,...,V_k){\rm d}V,\Re(\alpha_j)>{{p-1}\over2},\Re(\zeta_j)>{{p-1}\over2}.&(3.4)\cr}
$$Hence from (3.1) to (3.4) we can have the following theorem.

\vskip.3cm\noindent{\bf Theorem 3.1.}\hskip.3cm{\it Let $U_1,...,U_k$ be independently distributed as real matrix variate type-1 beta matrices with  parameters $(\zeta_j,\alpha_j),j=1,...,k, \Re(\alpha_j)>{{p-1}\over2},\Re(\zeta_j)>{{p-1}\over2}$. Let $V_1,..,V_k$ be another sequence of matrices having a joint density $f(V_1,...,V_k)$. Let the two sets $(X_1,...,X_k)$ and $(V_1,...,V_k)$ be independently distributed. Let $X_j=V_j^{1\over2}U_j^{-1}V_j^{1\over2},j=1,...,k$. Let $g(U_1,...,U_k)$ be the joint density of $U_1,...,U_k$. Then the Kober operator of the first kind for many matrix variables case as defined in (3.3) is a constant multiple of $g(U_1,...,U_k)$ as in (3.2).}

\vskip.3cm We can also have theorems parallel to the ones in section 2 and the proofs are parallel. Hence we list two such theorems here without proofs.

\vskip.3cm\noindent{\bf Theorem 3.2.}\hskip.3cm{\it Let $X_1,...,X_k$ have a joint real matrix-variate type-1 Dirichlet density with the parameters $(\zeta_1,...,\zeta_k;\zeta_{k+1})$. Consider the the transformation in (2.6) and let $Y_1,...,Y_k$ be defined there. Let $V_1,...,V_k$ be another sequence of matrix random variables having a joint density $f(V_1,...,V_k)$ where let $(X_1,...,X_k)$ and $(V_1,...,V_k)$ be independently distributed. Let $Y_j= V_j^{1\over2}U_j^{-1}V_j^{1\over2},j=1,...,k$. Let the joint density of $U_1,...,U_k$ be denoted by $g(U_1,...,U_k)$. Let

$$\gamma_j=\zeta_{j+1}+...+\zeta_{k+1}.\eqno(3.5)$$Then
$$\eqalignno{g(U_1,...,U_k)&=\{\prod_{j=1}^k{{\Gamma_p(\zeta_j+\gamma_j)}\over{\Gamma_p(\zeta_j)\Gamma_p(\gamma_j)}}\}
I_{U_1,...,U_k}^{(\zeta_j,\gamma_j),j=1,...,k}f(U_1,..,U_k)&(3.6)\cr
\noalign{\hbox{or}}
\{\prod_{j=1}^k{{\Gamma_p(\zeta_j)}\over{\Gamma_p(\zeta_j+\gamma_j)}}\}g(U_1,...,U_k)&=
I_{U_1,...,U_k}^{(\zeta_j,\gamma_j),j=1,...,k}f(U_1,...,U_k).&(3.7)\cr}
$$}

\vskip.3cm\noindent{\bf Theorem 3.3.}\hskip.3cm{\it Let $X_1,...X_k$ have a joint density as in (2.11) and the remaining transformations and notations remain as in Theorem 3.2. Let
$$\delta_j=\zeta_{j+1}+...+\zeta_{k+1}+\beta_j+...+\beta_k.\eqno(3.8)$$Let the joint density of $U_1,...,U_k$ be again denoted by $g(U_1,...,U_k)$. Then $g(U_1,...,U_k)$ is a density and
$$\{\prod_{j=1}^k{{\Gamma_p(\zeta_j)}\over{\Gamma_p(\zeta+\delta_j)}}\}g(U_1,...,U_k)
=I_{U_1,...,U_k}^{(\zeta_j,\delta_j),j=1,...,k}
f(U_1,...,U_k).\eqno(3.8)$$}

\vskip.3cm\noindent{\bf 4.\hskip.3cm M-transforms for the Fractional Integral Operators in the Many Matrix-variate Case}

\vskip.3cm Here we look at the M-transforms for fractional integral operators of the first and second kind in the many matrix-variate case. Consider the first kind operator in (3.4). The M-transform is given by
$$\eqalignno{M\{I_{U_1,...,U_K}^{(\zeta_j,\alpha_j),j=1,...k}&f(U_1,...,U_k)\}\cr
&=\int_{U_1>O}...\int_{U_k>O}|U_1|^{s_j-{{p+1}\over2}}...|U_k|^{s_k-{{p+1}\over2}}\cr
&\times I_{U_1,...,U_k}^{(\zeta_j,\alpha_j),j=1,...,k}f(U_1,...,U_k){\rm d}U_1\wedge...\wedge{\rm d}U_k\cr}
$$The $X$-integral is given by
$$\eqalignno{
\int_{U_j>V_j}|U_j|^{s_j-{{p+1}\over2}}&|U_j|^{-\zeta_j-\alpha_j}|U_j-V_j|^{\alpha_j-{{p+1}\over2}}{\rm d}U_j&=\int_{U_j>V_j}|V_j|^{-\zeta_j+s_j-{{p+1}\over2}}\int_{Y_j>O}|Y_j|^{\alpha_j-{{p+1}\over2}}\cr
&\times |I+Y_j|^{-(\zeta_j-s_j+{{p+1}\over2})}{\rm d}Y_j,~T_j=U_j-V_j,Y_j=V_j^{-{1\over2}}T_jV_j^{-{1\over2}}\cr
&=|V_j|^{-\zeta_j+s_j-{{p+1}\over2}}{{\Gamma_p(\alpha_j)\Gamma_p({{p+1}\over2}+\zeta_j-s_j)}
\over{\Gamma_p({{p+1}\over2}+\alpha_j+\zeta_j-s_j)}}\cr}
$$by evaluating the integral by using type-1 real matrix-variate beta integral, for $\Re(s)<\Re(\zeta+{{p+1}\over2}),\Re(\alpha_j)>{{p-1}\over2}$. Now the $V_j$-integrals give the M-transform of $f(V_1,...,V_k)$. Hence we can have the following theorem.

\vskip.3cm\noindent{\bf Theorem 4.1.}\hskip.3cm{\it For the fractional integral of the first kind of orders $\alpha_1,...,\alpha_k$ defined in (3.4) the M-transform is given by
$$\eqalignno{M\{I_{U_1,...,U_k}^{(\zeta_j,\alpha_j),j=1,...,k}f(U_1,...,U_k)\}&=f^{*}(s_1,...,s_k)\cr
&\times
\prod_{j=1}^k{{\Gamma_p({{p+1}\over2}+\zeta_j-s_j)}
\over{\Gamma_p({{p+1}\over2}+\alpha_j+\zeta_j-s_j)}}\cr}
$$for $\Re(s_j)<\Re(\zeta_j+{{p+1}\over2}),\Re(\alpha_j)>{{p-1}\over2},j=1,...,k.$, where $f^{*}(s_1,...,s_k)$ is the M=transform of $f(V_1,...,V_k)$.}

\vskip.3cm In a similar manner we can workout the M-transform of fractional integral operator of the second kind in the many matrix-variate case. In this context we start with (2.3). The M-tranform  is given by
$$\eqalignno{M\{K_{U_1,...,U_k}^{(\zeta_j,\alpha_j),j=1,...,k}f(U_1,...,U_k)\}
&=\{{{|U_j|^{\zeta_j}}\over{\Gamma_p(\alpha_j)}}\int_{U_j>O}|U_j|^{s_j-{{p+1}\over2}}\int_{V_j>U_j}
|V_j|^{-\zeta_j-\alpha_j}\cr
&\times |V_j-U_j|^{\alpha_j-{{p+1}\over2}}{\rm d}U_j\}f(V_1,...,V_k){\rm d}V_1\wedge...\wedge{\rm d}V_k.\cr}
$$The $U_j$-integral is given by
$$\eqalignno{\int_{U_j<V_j}|U_j|^{s_j-{{p+1}\over2}+\zeta_j}|V_j-U_j|^{\alpha_j-{{p+1}\over2}}{\rm d}U_j&=|V_j|^{\alpha_j-{{p+1}\over2}}\int_{U_j<V_j}|U_j|^{s_j+\zeta_j-{{p+1}\over2}}
|I-V^{-{1\over2}}U_jV_j^{-{1\over2}}|^{\alpha_j-{{p+1}\over2}}{\rm d}U_j.\cr
\noalign{\hbox{Put $Y_j=V^{-{1\over2}}U_jV_j^{-{1\over2}}$ and integrate out by using a matrix-variate type-1 beta integral}}
&=|V_j|^{\alpha_j-{{p+1}\over2}+s_j+\zeta_j}{{\Gamma_p(\alpha_j)\Gamma_p(\zeta_j+s_j)}
\over{\Gamma_p(\alpha_j+\zeta_j+s_j)}}.\cr}
$$Now the $v_j$-integrals give the M-transform of $f$. Hence we have the following theorem.

\vskip.3cm\noindent{\bf Theorem 5.2.}\hskip.3cm{\it For the fractional integral operator of the second kind defined in (2.3) the M-transform is given by
$$M\{K_{U_1,...,U_k}^{(\zeta_j,\alpha_j),j=1,...,k}f(U_1,...,U_k)\}=f^{*}(s_1,...,s_k)
\prod_{j=1}^k{{\Gamma_p(\zeta_j+s_j)}\over{\Gamma_p(\alpha_j+\zeta_j+s_j)}}
$$for $\Re(\alpha_j)>{{p-1}\over2},\Re(\zeta_j+s_j)>{{p-1}\over2},j=1,...,k$, where $f^{*}(s_1,...,s_k)$ is the M-transform of $f(V_1,...,V_k)$.}

\vskip.3cm Note that for $k=1$ the corresponding M-transforms in the one matrix variable case,  for $p=1$ the corresponding Mellin transforms in the $k$ scalar variables case and for $p=1,k=1$ the corresponding Mellin transforms in the one variable case for the Kober operators of the first  and second kinds are obtained from Theorems 5.1 and 5.2 respectively.

\vskip.3cm\noindent{\bf Acknowledgement}

\vskip.3cm The authors would like to thank the Department of Science and Technology, Government of India, New Delhi, for the financial assistance for this work under project number SR/S4/MS:287/05.

\vskip.3cm\centerline{\bf References}

\vskip.3cm\noindent Mathai, A.M. (1997):\hskip.3cm{\it Jacobians of Matrix Transformations and Functions of Matrix Argument}, World Scientific Publishing, New York.

\vskip.3cm\noindent Mathai, A.M. and Haubold, H.J. (2008):\hskip.3cm {\it Special Functions for Applied Scientists}, Springer, New York.

\bye